\documentclass[12pt]{article}
\usepackage{amsmath, amssymb, amsthm, euscript}

\theoremstyle{plain}
\newtheorem{lem}{Lemma}[section]
\newtheorem{thm}[lem]{Theorem} 
\newtheorem{cor}[lem]{Corollary}
\newtheorem{prop}[lem]{Proposition}

\theoremstyle{remark}
\newtheorem{rem}[lem]{Remark}

\theoremstyle{definition}
\newtheorem{exm}[lem]{Example}
\newtheorem{defi}[lem]{Definition}
\newtheorem{para}[lem]{}

\newcommand{\vp}{\varepsilon}

\newcommand{\n}{\noindent}
\newcommand{\cl}[1]{{\mathcal{#1}}}
\newcommand{\bb}[1]{{\mathbb{#1}}}
\newcommand\Ac{{\mathcal{A}}}
\newcommand\Bc{{\mathcal{B}}}
\newcommand\Cpx{{\mathbb C}}
\newcommand\Fb{{\mathbb F}}
\newcommand\HEu{{\EuScript H}}                   
\newcommand\HEuo{\HEu^\circ}                   
\newcommand\id{{\operatorname{id}}}
\newcommand\Ints{{\mathbb Z}}
\newcommand\KEu{{\EuScript K}}                   
\newcommand\Nats{{\mathbb N}}
\newcommand\oneh{{\hat 1}}
\newcommand\otdt{\otimes\cdots\otimes}
\newcommand\Puk{{\operatorname{Puk}}}
\newcommand\tauhat{{\hat\tau}}
\newcommand\Vc{{\mathcal{V}}}
\newcommand\pit{{\tilde\pi}}
\newcommand\mut{{\tilde\mu}}
\newcommand\mt{{\tilde m}}
\newcommand\Hh{{\widehat H}}
\newcommand\oup{^{\mathrm o}}
\numberwithin{equation}{section}
\newcommand\xpt{\widetilde{x'}}
\newcommand\xt{{\tilde x}}
\newcommand\yt{{\tilde y}}
\newcommand\ut{{\tilde u}}
\newcommand\Bbar{{\overline B}}
\newcommand\Abar{{\overline A}}
\newcommand\ah{{\hat a}}
\newcommand\Afr{{\mathfrak A}}
\newcommand\MM{measure--multiplicity }
\begin{document}

\title{VALUES OF THE PUK\'ANSZKY INVARIANT\\ IN FREE GROUP FACTORS AND\\ THE
HYPERFINITE FACTOR}

\author{}
\date{}
\maketitle

\begin{center}
\begin{tabular}{lll}
KENNETH J.\ DYKEMA$^{(*)}$&ALLAN M.\ SINCLAIR&ROGER R.\ SMITH$^{(**)}$\\
Department of Mathematics&School of Mathematics&Department of 
Mathematics\\
Texas A\&M University&University of Edinburgh&Texas A\&M University\\
College Station, TX 77843&Edinburgh, EH9 3JZ&College Station, TX 77843\\
USA&SCOTLAND&USA\\
{\tt kdykema@math.tamu.edu}&{\tt A.Sinclair@ed.ac.uk}&{\tt 
rsmith@math.tamu.edu}
\end{tabular}
\end{center}

\begin{abstract}
Let $A\subseteq M \subseteq B(L^2(M))$ be a maximal abelian self-adjoint
subalgebra
(masa) in a type ${\text{II}}_1$ factor $M$ in its standard representation.
The abelian von~Neumann algebra $\cl A$ generated by $A$ and $JAJ$ has a type
$\text{I}$ commutant which contains the projection $e_A\in \cl A$ onto
$L^2(A)$.
Then $\cl A'(1-e_A)$ decomposes into a direct sum of type ${\text{I}}_n$
algebras for $n\in \{1,2,\cdots,\infty\}$, and those $n$'s which occur in the
direct sum form a set called the Puk\'anszky invariant, $\Puk(A)$, also
denoted $\Puk_M(A)$ when the containing factor is ambiguous. In this paper we
show that this invariant can take on the values $S\cup \{\infty\}$ when $M$ is
both a free group factor and the hyperfinite factor, and where $S$ is an
arbitrary subset of $\Nats$. The only
previously known values for masas in free group factors were $\{\infty\}$ and
$\{1,\infty\}$,
and some values of the form $S\cup\{\infty\}$ are new also for the hyperfinite
factor.

We also consider a more refined invariant
(that we will call the \MM invariant),
which was considered
recently by Neshveyev and St\o{}rmer and has been known to experts for a long
time.
We use the \MM invariant to distinguish two masas in a free group factor,
both having Puk\'anszky invariant $\{n,\infty\}$, for arbitrary $n\in\Nats$.
\end{abstract}

\vfill

\n $\underline{\hspace{1in}}$

\n {\footnotesize $^{(*)}$Research supported in part by NSF grant DMS-0300336, 
 the Alexander von Humboldt Foundation, and the Edinburgh Mathematical Society.

\n $^{(**)}$Research supported in part by NSF grant DMS-0401043.}
\newpage

\section{Introduction}\label{sec1}

\indent

The Puk\'{a}nszky invariant ${\text{Puk}}(A)$ of a maximal abelian
self-adjoint
subalgebra (masa) $A$ of a separable type ${\text{II}}_1$ factor $N$ with
normalized 
normal trace $\tau$ was introduced in \cite{Puk}. If $J$ denotes the canonical
involution on $L^2(N,\tau)$, then the abelian von Neumann algebra $\cl A$
generated by $A$ and $JAJ$ has a type ${\text{I}}$ commutant $\cl A'$. The
projection $e_A$ onto $L^2(A)$ lies in $\cl A$ and $\cl A'(1-e_A)$
decomposes into a direct sum of type ${\text{I}}_{n_i}$ algebras, where $1\leq
n_i\leq \infty$. Those $n_i$'s in this sum form  ${\text{Puk}}(A)$, a subset
of $\bb N \cup \{\infty\}$. This quantity is invariant under the action of any
automorphism of $N$, and so serves as an aid to distinguishing pairs of masas.
In
\cite{Puk,NS,SSpuk}, various values of the invariant were found for masas,
primarily
in the hyperfinite type ${\text{II}}_1$ factor $R$. In this paper we consider
the possible values of the invariant for masas in the free group factors
$L(\bb F_n)$ for $2\leq n \leq \infty$.  In
this paper one of our main objectives is to show that
in free group factors
there exist strongly singular masas
whose Puk\'{a}nszky invariants are $S\cup \{\infty\}$ where
$S$  is an arbitrary subset of $\bb N$ (this and other terminology will be
explained in the next
section).
There are two standard examples of masas in the free group factors. One arises
from a single generator of $\bb F_k$ and the criteria of \cite{Dix} show that
it is singular, while its invariant is $\{\infty\}$, \cite{Rob}. 
The other type is the radial or laplacian masa. R\u{a}dulescu, \cite{Rad},
has shown that this masa in $L(\bb F_2)$ has $\{\infty\}$
for its invariant, and it was shown in \cite{SSgafa,SStrans} that both types of
masas are
strongly singular. On the other hand, the isomorphism between $L(\bb F_2)$ and
$L(\bb
F_5)\otimes \bb M_2$, \cite{Voi1}, can be used to find a tensor product masa
whose invariant is $\{1,\infty\}$, although this masa is not singular. These
two possibilities, $\{\infty\}$ and $\{1,\infty\}$, were previously the only
ones known,
although it was shown in \cite{SSpuk} that the invariant in a free group
factor cannot be a finite set of integers, based on results from \cite{D2}.

In the second section we investigate masas in free product algebras $M*Q$ where
$M$ is a diffuse finite von Neumann algebra while $Q$ has a finite trace and
dimension at least two. We show that any (singular) masa in $M$ is also a
(singular) masa in $M*Q$ by determining equality of the unitary normalizers of
such a masa in the two algebras. This result was originally obtained in
\cite{Popa} by different methods. We then recall from \cite{SSpuk} the 
construction of masas
$A_n$, $1\leq n <\infty$, with ${\text{Puk}}(A_n)=\{n,\infty\}$, and show that
they
are strongly singular in their containing factors $M_n$ and in $M_n*Q$.    
These masas form the basis for the examples of the next section.

In the third section, for each set
$S\subseteq \bb N$, we form the masa $A$ as a direct sum of the masas $\{A_i\}_
{i\in S}$ inside the direct sum of $\{M_i\}_{i\in S}$. The free product with a
suitably chosen von Neumann algebra $Q$ gives a masa in $L(\bb F_2)$ and the
main results of the section, Theorems \ref{thm3.2} and \ref{thm3.3}, determine
the Puk\'{a}nszky
invariant as $S\cup \{\infty\}$. A further free product with $L(\bb F_{n-2})$
then gives the corresponding result for each $L(\bb F_n)$, $3\leq n \leq
\infty$. More generally, we show that the same conclusions also hold for $L(\bb
F_n*\Gamma)$, where $\Gamma$ is an arbitrary countable discrete group.

In the
fourth section we obtain the same existence result
for masas having Puk\'{a}nszky invariant $S\cup\{\infty\}$, but
in the hyperfinite factor.
This extends the set of known values from those obtained in \cite{NS,SSpuk}.

In the last
two sections we consider a finer invariant which can distinguish masas with
the same Puk\'{a}nszky invariant. This has been used extensively in \cite{NS},
where its origins are attributed to \cite{FM}. Since it involves a measure
class and a multiplicity function, we will refer to it as the
\MM invariant. We use it to give examples in the sixth
section of nonconjugate masas in a free group factor,
both with Puk\'{a}nszky invariant $\{n,\infty\}$.
The groundwork for this result is laid in the fifth section.

We assume throughout, even when not stated explicitly, that all von~Neumann
algebras are separable in the sense that their preduals are norm-separable as
Banach spaces. This is equivalent to the assumption of faithful representations
on separable Hilbert spaces. Moreover, all groups will be discrete and
countable, so that the associated group von~Neumann algebras will be separable.

The results of this paper rely heavily on the theory of free products of von
Neumann algebras. We refer the reader to \cite{VDN} for the necessary
background material, and also to \cite{BR,NT,Rad} for related results on masas
in free group factors. We note that {\em singularity} and {\em strong
singularity} for masas were shown recently to be equivalent, \cite{SSWW}. We
have retained the latter terminology in this paper, since showing strong
singularity is often the most direct way of proving singularity of masas.

\medskip

\noindent{\bf Acknowledgements.}
Much of this work was accomplished while KJD was visiting   
the University of M\"unster during 2004-2005. The warm
hospitality of the members of the Mathematics Institute is gratefully
acknowledged.
KJD and RRS thank Kunal Mukherjee for helpful comments about early versions of
this paper.


\section{Masas in free product factors}\label{sect2}

\indent

Throughout this section $M$ and $Q$ will be  separable von Neumann algebras
with 
faithful normal tracial states $\tau_M$ and $\tau_Q$. We assume that $M$ is 
diffuse and that $Q$ has dimension at least 2. Moreover, $M$ and $Q$ are in 
standard form on separable Hilbert spaces $\HEu_1$ and $\HEu_2$ respectively so 
that there are distinguished unit vectors $\xi_i\in \HEu_i$, $i=1,2$, such that 
$\tau_M(\cdot) = \langle\cdot\xi_1,\xi_1\rangle$ and $\tau_Q(\cdot) = \langle 
\cdot\xi_2,\xi_2\rangle$. Recall from \cite{VDN} that the free product
von~Neumann algebra $M*Q$ has a normal tracial state $\tau_{M*Q}$ whose
restrictions to $M$ and $Q$ are respectively $\tau_M$ and $\tau_Q$. 

If $A$ is a diffuse abelian von~Neumann subalgebra of $M$ then it is isomorphic 
to $L^\infty[0,1]$ with a faithful state induced from the trace on $M$. This 
corresponds to a probability measure $\mu$ of the form $f(t)\,dt$ where $f\in 
L^1[0,1]^+$ has unit norm. Then $(L^\infty[0,1], dt)$ and $(L^\infty[0,1], f\,
dt)$ are continuous masas on $(L^2[0,1], dt)$ and $(L^2[0,1], f \,dt)$ 
respectively, and they are unitarily equivalent by a unitary which takes one 
separating and cyclic vector 1 to the other $f^{1/2}$, 
\cite[Theorem 9.4.1]{KR}. The 
unitaries $\{e^{int}\}_{n\in {\bb Z}}$ then give rise to a set of unitaries 
$\{v^n\}_{n\in {\bb Z}}$  satisfying the orthogonality condition
\begin{equation}\label{eq2.1}
\tau_M(v^{*m}v^n) = \tau_M(v^{n-m}) = 0,\qquad m\ne n.
\end{equation}
Such an operator $v$ is called a {\it{Haar unitary}}.

Recall from \cite{SSgafa} that for a map $\phi\colon \ N\to N$, where $N$ is 
a von~Neumann algebra with a faithful normal trace $\tau$, 
$\|\phi\|_{\infty,2}$ is defined by
\begin{equation}\label{eq2.2}
\|\phi\|_{\infty,2} = \sup\{\|\phi(x)\|_{2,\tau} \colon \ \|x\|\le 1, 
\ x\in 
N\}.
\end{equation}
The subscript $\tau$ indicates the trace used to calculate $\|\cdot\|_2$ when
there could be ambiguity.
We will require maps of the form $\Phi_{s,t}(x) = sxt$, $x\in N$, where $s$ 
and $t$ are fixed but arbitrary elements of $N$. The following simple lemma 
will be useful subsequently, and indicates that $\Phi_{s,t}$ depends 
continuously on $s$ and $t$ in an appropriate sense.

\begin{lem}\label{lem2.1}
If $s,s',t,t'\in N$ are operators of norm at most 1, then
\begin{equation}\label{eq2.2a}
\|\Phi_{s,s'}-\Phi_{t,t'}\|_{\infty,2} \le \|s-t\|_2+ \|s'-t'\|_2. 
\end{equation}
\end{lem}

\begin{proof}
This is immediate from the algebraic identity

\begin{equation}\label{eq2.3}
sxs'-txt'=(s-t)xs'+tx(s'-t'),
\end{equation}
for $x\in N$.
\end{proof}

\begin{lem}\label{lem2.2}
Let ${\bb E}^{M*Q}_M$ be 
the unique trace preserving conditional expectation of $M*Q$ onto $M$ and let 
$v$ be a Haar unitary in $M$. Then
\begin{equation}
\label{eq2.4}
\lim_{|k|\to\infty} \|{\bb E}^{M*Q}_M(xv^ky) - {\bb E}^{M*Q}_M(x) v^k {\bb 
E}^{M*Q}_M(y)\|_{2,\tau_M} = 0
\end{equation}
for all $x,y\in M*Q$.
\end{lem}

\begin{proof}
In view of Lemma \ref{lem2.1}, it suffices to prove the result when $x$ and 
$y$ are words in $M*Q$. Each such word $w$ can be expressed as $w=m + 
\sum\limits^r_{i=1} \alpha_iw_i$ where $m\in M$, $\alpha_i\in {\bb C}$, and 
each $w_i$ contains only letters of zero trace, and at least one letter from
$Q$, \cite[Lemma 1]{Chi}. For such a
representation, ${\bb E}^{M*Q}_M(w)=m$.
Now \eqref{eq2.4} 
is immediate if $x$ or $y$ lies in $M$, by properties of the conditional 
expectation. Thus, it suffices to prove \eqref{eq2.4} when $x$ and $y$ are 
words whose letters have zero trace and contain at least one letter from $Q$. 
There are several cases to consider, depending on whether the last letter of 
$x$ and the first letter of $y$ are in $M$ or $Q$. Suppose initially that 
both are in $Q$. Then $xv^ky$ is a reduced word for $|k|\ge 1$, each letter 
has zero trace, and thus the two terms in \eqref{eq2.4} are both 0. Now suppose 
that $x$ ends and $y$ begins with letters from $M$. Then write $x=\tilde 
xm_1$, $y = m_2\tilde y$ where $m_1,m_2\in M$, and the last letter of $\tilde 
x$ and the first letter of $\tilde y$ are in $Q$. Then we may express
\begin{equation}\label{eq2.5}
xv^ky = \tilde xm_1v^k m_2\tilde y = \tilde x(m_1v^km_2-\tau_M(m_1v^km_2)) 
\tilde y + \tau_M(m_1v^km_2) \tilde x\tilde y.
\end{equation}
Thus
\begin{equation}\label{eq2.6}
{\bb E}^{M*Q}_M(x) = {\bb E}^{M*Q}_M(y) = {\bb E}^{M*Q}_M(\tilde x(m_1v^k 
m_2-\tau_M(m_1v^km_2))\tilde y) = 0,
\end{equation}
while $\lim\limits_{|k|\to\infty} \tau_M(m_2m_1v^k)=0$ since $\{v^k\}_{k\in 
{\bb Z}}$ is an orthonormal set of vectors in $L^2(M,\tau_M)$. This proves 
\eqref{eq2.4}  in this case also. The remaining case, when exactly one 
of $x$ and $y^*$ ends with a letter from $Q$, is handled just as in the 
previous case but using only one of $\tilde x$ or $\tilde y$.
\end{proof}

We denote by ${\cl N}_M(A)$ the set of unitaries $u\in M$ which normalize a 
given subalgebra $A$, in the sense that $uAu^*=A$. This group  is called the
{\it{unitary normalizer}}
of $A$ in $M$. The masa $A$ is said to be
{\it{singular}} if ${\cl N}_M(A)$ coincides
with the unitary group of $A$, \cite{Dix}.

The following result is due to Popa, \cite[Remark 6.3]{Popa}. We present an
alternative proof, based on conditional expectations, since this method is
required below.

\begin{thm}\label{thm2.3}
Let $A$ be a diffuse von~Neumann subalgebra of $M$. Then the following
statements hold:
\begin{itemize}
\item[(i)] The unitary normalizers ${\cl N}_M(A)$ and ${\cl 
N}_{M*Q}(A)$ are equal. 
\item[(ii)] If $A$ is a (singular) masa in $M$, then it is also a 
(singular) masa in $M*Q$.
\end{itemize}
\end{thm}

\begin{proof}
(i) Let $u\in M*Q$ be a unitary which normalizes $A$, and let $v$ be a Haar 
unitary in $A$. Then $uv^ku^*\in A\subseteq M$ for $k\in {\bb Z}$, and so 
${\bb E}^{M*Q}_M(uv^ku^*) = uv^ku^*$. Thus the $\|\cdot\|_{2,\tau_M}$-norms 
of these elements are 1. From \eqref{eq2.4} in Lemma \ref{lem2.2},
\begin{equation}\label{eq2.7}
\lim_{|k|\to\infty} \|{\bb E}^{M*Q}_M(u)v^k {\bb E}^{M*Q}_M(u^*)\|_{2,\tau_M} 
= 1,
\end{equation}
which is impossible unless $\|{\bb E}^{M*Q}_M(u)\|_{2,\tau_M} = 1$. But then 
Hilbert space orthogonality gives
\begin{equation}\label{eq2.8}
1 = \|u\|^2_{2,\tau_{M*Q}} = \|(I-{\bb E}^{M*Q}_M)(u)\|^2_{2,\tau_{M*Q}} + 
\|{\bb E}^{M*Q}_M(u)\|^2_{2,\tau_M},
\end{equation}
from which we conclude that $(I-{\bb E}^{M*Q}_M)(u) = 0$. Thus $u\in M$, 
proving that ${\cl N}_{M*Q}(A) \subseteq {\cl N}_M(A)$. The reverse 
containment is obvious.

(ii) Any unitary  $u\in M*Q$ which commutes with $A$ lies in ${\cl
N}_{M*Q}(A)$, and 
so in $M$, by part (i). If $A$ is a masa in $M$, then $u\in A'\cap M = A$, and
so $A$ is a 
masa in $M*Q$. Further, suppose that $A$ is singular. Then each unitary in 
${\cl N}_M(A)$ lies in $A$ and so the same is true for ${\cl N}_{M*Q}(A)$. 
Thus $A$ is also singular in $M*Q$.
\end{proof}

In \cite{SSgafa}, {\it{strong singularity}} of a masa $A\subseteq
M$ 
was defined by the requirement that
\begin{equation}\label{eq2.9}
\|{\bb E}^M_{uAu^*} - {\bb E}^M_A\|_{\infty,2} \ge \|u-{\bb 
E}^M_A(u)\|_{2}
\end{equation}
for all unitaries $u\in M$. The left--hand side of \eqref{eq2.9} vanishes when 
$u\in {\cl N}_M(A)$, and the inequality then shows that $u\in A$. Thus  
singularity of $A$ is a consequence of strong singularity, and the reverse
implication was established recently in \cite{SSWW}.  The usefulness of 
strong singularity lies in the ease with which \eqref{eq2.9} can be verified 
in specific cases. There are two main criteria which establish strong
singularity. If there is a unitary $v\in A$ such that
\begin{equation}\label{eq2.10}
\lim_{|k|\to\infty} \|{\bb E}^M_A(xv^ky) - {\bb E}^M_A(x) v^k{\bb 
E}^M_A(y)\|_2 = 0
\end{equation}
for all $x,y\in M$, then we say that ${\bb E}^M_A$ is an {\it{asymptotic 
homomorphism}} with respect to $v$. In \cite[Theorem 4.7]{SSgafa} it was shown 
that this property implies strong singularity for $A$ when $M$ is a type 
$\text{II}_1$ factor, but the proof is also valid for a general finite 
von~Neumann algebra. There is a weaker form of \eqref{eq2.10}, defining the
{{\it weak asymptotic homomorphism property}} ({\it{WAHP}}): given $\vp >0$
and a finite set of elements $\{x_1,\ldots ,x_n,y_1,\ldots ,y_n\}\subseteq M$,
there exists a unitary $u\in A$ such that
\begin{equation}\label{eq2.10z}
\|{\bb E}^M_A(x_iuy_j) - {\bb E}^M_A(x_i) u{\bb 
E}^M_A(y_j)\|_2 <\vp,\ \ \ 1\leq i,j\leq n.
\end{equation}
This property implies strong singularity, \cite[Lemma 2.1]{RSS}, and it is also
equivalent to singularity, \cite[Theorem 2.3]{SSWW}. The previous theorem and
these remarks lead to the following, which we include just to emphasize the
methods employed for subsequent results.
\begin{cor}\label{cor2.4}
Let $A$ be a  masa in $M$. 
\begin{itemize}
\item[(i)] If ${\bb E}^M_A$ is an 
asymptotic homomorphism with respect  to a Haar unitary $v\in A$ then this also
holds for ${\bb E}^{M*Q}_A$;
\item[(ii)] if $A$ has the {\it WAHP} in $M$ then it also has this property in
$M*Q$.
\end{itemize}
In both cases $A$ is strongly singular in both $M$ and $M*Q$.
\end{cor}

\begin{proof}
(i) We first note that uniqueness of trace preserving conditional expectations 
implies that ${\bb E}^{M*Q}_A = {\bb E}^M_A \circ {\bb E}^{M*Q}_M$. From 
Lemma \ref{lem2.2},
\begin{equation}\label{eq2.11}
\lim_{|k|\to\infty} \|{\bb E}^{M*Q}_M(xv^ky) - {\bb E}^{M*Q}_M(x) v^k{\bb 
E}^{M*Q}_M(y)\|_2  = 0
\end{equation}
for all $x,y\in M*Q$. Since ${\bb E}^M_A$ is a $\|\cdot\|_2$-norm 
contraction, we may apply this operator to \eqref{eq2.11} to obtain
\begin{equation}\label{eq2.12}
\lim_{|k|\to\infty} \|{\bb E}^{M*Q}_A(xv^ky) - {\bb E}^M_A({\bb E}^{M*Q}_M(x) 
v^k{\bb E}^{M*Q}_M(y))\|_2 = 0
\end{equation}
for all $x,y\in M*Q$. The asymptotic homomorphism hypothesis, when applied to 
the second term in \eqref{eq2.12}, leads to
\begin{equation}\label{eq2.13}
\lim_{|k|\to\infty} \|{\bb E}^{M*Q}_A(xv^ky) - {\bb E}^{M*Q}_A(x)v^k {\bb 
E}^{M*Q}_A(y)\|_2 =  0
\end{equation}
for all $x,y\in M*Q$. Thus ${\bb E}^{M*Q}_A$ is an asymptotic homomorphism 
with respect to $v$. 

\noindent (ii) This is immediate from Theorem \ref{thm2.3} and the equivalence
of singularity and the {\it WAHP}.

 In both cases, the (strong) singularity of $A$ follows from the remarks
preceding this corollary.
\end{proof}

In \cite{SSpuk}, a family of masas inside group von Neumann factors was 
presented with various Puk\'anszky invariants. We recall these masas now 
since we wish to give some extra information about them. Consider an abelian
subgroup $H$ of an I.C.C. group $G$. In \cite{SSpuk}, the problem of
describing the Puk\'anszky invariant of $L(H)$ in $L(G)$, when $L(H)$ is a
masa, was solved in terms of double cosets and stabilizer subgroups. The
double coset $HgH$ of $g\in G\backslash H$ is $\{hgk\colon\ h,k\in H\}$. The
stabilizer subgroup $K_g\subseteq H\times H$ is defined to be $\{(h,k)\colon\
h,k\in H,\
hgk=g\}$. Under the additional hypothesis that any two such subgroups $K_c$ and
$K_d$, for $c,d\in G\backslash H$, are either equal or satisfy the
noncommensurability condition that
$K_cK_d/(K_c\cap K_d)$ has infinite order, an equivalence relation was defined
on the nontrivial double cosets by $HcH\sim HdH$ if and only if $K_c=K_d$. The
numbers, including $\infty$, in ${\text{Puk}}(L(H))$ are then the numbers of
double cosets that occur in the various equivalence classes,
\cite[Theorem 4.1]{SSpuk}. We use this now to discuss
the Puk\'anszky invariant in certain free products of group factors. The
groups that arise all satisfy the additional hypothesis above, so that
\cite[Theorem 4.1]{SSpuk} applies to them.

\begin{prop}\label{prop2.5}
Let $\Gamma$ be a countable discrete group of order at least 2. For each $n\geq
1$, there exists a countable discrete amenable i.c.c.\ group $G_n$ with an
abelian
subgroup $H_n$ having the following properties:
\begin{itemize}
\item[(i)] $L(H_n)$ is a strongly singular masa in both $L(G_n)$ and
$L(G_n)*L(\Gamma)$.
\item[(ii)]
$\Puk_{L(G_n)}(L(H_n))=\Puk_{L(G_n)*L(\Gamma)}(L(H_n))=\{n,\infty\}$.
\end{itemize}
\end{prop}
\begin{proof}
Fix an integer $n\in {\bb N}$. Let ${\bb Q}$ denote the group of rationals 
under addition and let ${\bb Q}^\times$ be the multiplicative group of 
nonzero rationals. Let $P_n\subseteq {\bb Q}^\times$ be the subgroup
\[
P_n = \left\{\frac{p}q 2^{nj}\colon \ j\in {\bb Z},\  p,q\in {\bb 
Z}_{\text{odd}}\right\},\qquad 1\leq n < \infty,
\]
let $P_{\infty}\subseteq {\bb Q}^\times$ be the subgroup
\[
P_{\infty} = \left\{\frac{p}q \colon \  p,q\in {\bb 
Z}_{\text{odd}}\right\},
\]
and let $G_n$, $n\geq 1$, be the matrix group
\[G_n = \left\{\left(\begin{matrix} 1&x&y\\ 0&f&0\\
0&0&g\end{matrix}\right)\colon 
\ 
x,y\in {\bb Q}, \ f\in P_n, \ g\in P_{\infty}\right\}
\]
with abelian subgroup $H_n$ consisting of the diagonal matrices in $G_n$. Then 
$L(H_n)$ is a masa in the factor $ L(G_n)$, and $\text{Puk}(L(H_n )) 
= \{n,\infty\}$, (see \cite[Example 5.2]{SSpuk}, where it was also noted that
$G_n$ is amenable). There are three equivalence classes
of double cosets whose sizes are $n$, $\infty$ and $\infty$, and the
corresponding stabilizer subgroups are respectively
\begin{align}
&\left\{\left(\left(\begin{matrix} 1&0&0\\ 0&1&0\\ 0&0&g\end{matrix}\right), 
\left(\begin{matrix} 1&0&0\\ 0&1&0\\ 0&0&g^{-1}\end{matrix}\right)\right)\colon 
\ g\in P_{\infty}\right\}, \left\{\left(\left(\begin{matrix} 1&0&0\\ 0&f&0\\ 
0&0&1\end{matrix}\right), \left(\begin{matrix} 1&0&0\\ 0&f^{-1}&0\\ 
0&0&1\end{matrix}\right)\right)\colon \ f\in P_n\right\},\nonumber\\
\label{eq5.16}
&\qquad \left\{\left(\left(\begin{matrix} 1&0&0\\ 0&1&0\\ 
0&0&1\end{matrix}\right), \left(\begin{matrix} 1&0&0\\ 0&1&0\\ 
0&0&1\end{matrix}\right)\right)\right\}.
\end{align}
When $H_n$ is viewed as a subgroup of $G_n*\Gamma$, where $\Gamma$ is a
countable discrete group of order at least 2, each element of
$(G_n*\Gamma)\backslash G_n$ has a trivial stabilizer subgroup so the extra
double
cosets fall into the third equivalence class above, showing that
\[{\text{Puk}}_{L(G_n)*L(\Gamma)}(L(H_n))={\text{Puk}}_{L(G_n)}(L(H_n))=\{n,\infty\}\]
in this case.
Let $v_n$ denote the group element 
$\left(\begin{smallmatrix} 1&0&0\\ 0&3&0\\0&0&5\end{smallmatrix}\right)\in
H_n$,
viewed as a 
unitary in $L(H_n)$. If $\tau_n$ is the standard trace on $L(G_n)$, then 
$\tau_n(g) = 0$ for any group element $g\ne e$, so $v_n$ is a Haar unitary in 
$L(H_n)$. 
A routine matrix calculation shows that  $g_1v^k_ng_2 \in G_n\backslash H_n$
for $|k|$ sufficiently large, when $g_1,g_2\in G_n\backslash H_n$.  
Then, viewing these group elements as unitaries in $L(G_n)$, we have ${\bb 
E}^{L(G_n)}_{L(H_n)}(g_i) = 0$, $i=1,2$, and ${\bb
E}^{L(G_n)}_{L(H_n)}(g_1v^k_ng_2)  = 0$ 
for $|k|$ sufficiently large, so that
\begin{equation}\label{eq2.15}
\lim_{|k|\to\infty} \|{\bb E}^{L(G_n)}_{L(H_n)}(g_1v^k_ng_2) - {\bb 
E}^{L(G_n)}_{L(H_n)}(g_1) v^k_n{\bb E}^{L(G_n)}_{L(H_n)}(g_2)\|_2 = 0.
\end{equation}
Then a  simple approximation argument gives \eqref{eq2.15} for $g_1$ and 
$g_2$ replaced by general elements $x,y\in L(G_n)$. This shows that ${\bb 
E}^{L(G_n)}_{L(H_n)}$ is an asymptotic homomorphism for $v_n$, where we have 
employed the methods of \cite[Section 2]{RSS}. Thus each 
$L(H_n)$ is strongly singular in both $L(G_n)$ and $L(G_n)*L(\Gamma)$, using
Corollary 
\ref{cor2.4}.
\end{proof}

\begin{rem}\label{rem2.5a}
The groups constructed above will form the basis for our results in the next
section.
Denote the von~Neumann algebras $L(G_n)$ and $L(H_n)$ of Proposition
\ref{prop2.5} by $M_n$ and $A_n$ respectively.
Now consider an arbitrary nonempty subset $S\subseteq {\bb N}$, let $M_S = 
\bigoplus\limits_{n\in S} M_n$ and let $A_S = \bigoplus\limits_{n\in S} A_n$. 
Choose numbers $\{\alpha_n\}_{n\in S}$ from (0,1) such that 
$\sum\limits_{n\in S} \alpha_n=1$, and let $\tau_S$ be the normalized trace 
on $M_S$ given by $\tau_S = \sum\limits_{n\in S}\alpha_n\tau_n$, where 
$\tau_n$ is the normalized trace on $M_n$, $n\in S$. Then let $v_S = 
\bigoplus\limits_{n\in S} v_n$ be a unitary in $M_S$. Since arbitrary direct 
sums are approximable in $\|\cdot\|_{2,\tau_S}$-norm by finitely nonzero 
ones, it is routine to verify that ${\bb E}^{M_S}_{A_S}$ is an asymptotic 
homomorphism for $v_S$ and so, by Corollary \ref{cor2.4}, $A_S$ is a strongly 
singular masa in $M_S*Q$ for any von Neumann algebra $Q$ of dimension at least
2. These examples will be used below with $Q=L(\Gamma)$ for various choices of
$\Gamma$.$\hfill\square$
\end{rem}

We end this section with two observations which we include because they are 
simple deductions from our previous work. The first is known, \cite{D,D1}, and
the 
second may be known but we do not have a reference.

\begin{rem}\label{rem2.6}
(i) If $M$ is a diffuse finite von Neumann algebra then any central unitary
$u\in 
Z(M*Q)$ normalizes all masas in $M$ so, by Theorem \ref{thm2.3}, $u\in M$. 
From properties of the free product, the only elements in $M$ which commute 
with $Q$ are scalar multiples of 1, so $Z(M*Q)$ is trivial and $M*Q$ is a 
factor.

(ii) With the same assumption on $M$, and $Q$ any type $\text{II}_1$ factor,
$M$ 
embeds into the factor $M*Q$ with $\tau_M$ being the restriction to $M$ of 
$\tau_{M*Q}$. Moreover, if $u$ is a unitary in $M'\cap (M*Q)$, then $u$ 
normalizes each masa in $M$, so must lie in $M$, by Theorem \ref{thm2.3}. It 
follows that, for this embedding, the relative commutant and the center 
coincide. Of course, $Z(M)$ will be present in the relative commutant for  
any embedding of $M$ into a factor.$\hfill\square$
\end{rem}

\section{Construction of masas}

\indent

 In this section we construct strongly singular masas in the free group
factors whose Puk\'anszky invariants are $S\cup\{\infty\}$ for arbitrary
subsets $S$
of $\Nats$. The construction is based on direct sums, and the following two
results keep track of the contributions of the individual summands.

\begin{lem}\label{lem3.1}
Suppose that $B$ and $D$ are separable von Neumann algebras with normal
faithful traces
$\tau_B$ and $\tau_D$, respectively.
Suppose that $B$ and $D$ are both of dimension at least two.
Let $\lambda_B$ and $\rho_B$ denote the usual left and right actions of $B$ on
$L^2(B):=L^2(B,\tau_B)$.
Let
\[
(N,\tau)=(B,\tau_B)*(D,\tau_D)
\]
be their free product von Neumann algebra.
Let $\lambda_N$ and $\rho_N$ denote as usual the left and right actions of $N$
on
$L^2(N):=L^2(N,\tau)$.
Then there is a separable infinite dimensional Hilbert space $\KEu$ and a
unitary operator
\begin{equation}\label{eq:U}
U:L^2(N)\to L^2(B)\;\oplus\; L^2(B)\otimes\KEu\otimes L^2(B)
\end{equation}
such that for all $b\in B$,
\begin{align}
U\lambda_N(b)U^*&=\lambda_B(b)\;\oplus\;
\lambda_B(b)\otimes\id_{\KEu}\otimes\id_{L^2(B)} \label{eq:Ulam} \\
U\rho_N(b)U^*&=\rho_B(b)\;\oplus\;
\id_{L^2(B)}\otimes\id_{\KEu}\otimes\rho_B(b). \label{eq:Urho}
\end{align}
\end{lem}
\begin{proof}
This is very similar to part of Voiculescu's construction of the free product
of von Neumann
algebras~\cite{V85}, but for completeness we will
describe it in some detail.
Let us write
$\HEu_D=L^2(D,\tau_D)$ and let $\HEuo_D$ denote the orthocomplement of the
specified
vector $\oneh_D\in\HEu_D$.
Similarly, let $\HEu_B=L^2(B)$ and let $\HEuo_B$ denote the orthocomplement of
the specified
vector $\oneh_B\in\HEu_B$.
Then by Voiculescu's construction,
\begin{equation}\label{eq:L2M}
L^2(N)=\Cpx\xi\oplus
\bigoplus_{\substack{n\in\Nats \\ i_1,\ldots,i_n\in\{B,D\} \\
i_1\ne i_2,\,i_2\ne i_3,\,\ldots,\,i_{n-1}\ne i_n}}
\HEuo_{i_1}\otimes\cdots\otimes\HEuo_{i_n}.
\end{equation}
Let
\[
\KEu=\HEuo_D\oplus\bigoplus_{n=1}^\infty\HEuo_D\otimes(\HEuo_B\otimes\HEuo_D)^{\otimes
n}
\]
and consider the unitary
\[
U:L^2(N)\to L^2(B)\;\oplus\; L^2(B)\otimes\KEu\otimes L^2(B)
\]
defined, using~\eqref{eq:L2M}, as follows. The distinguished vector $\xi$ is
mapped to $\oneh_B\in L^2(B)$ and $U$ acts as the identity on $\HEuo_B
\subseteq L^2(B)$. For   vectors $\zeta_1\otdt\zeta_n \in
\HEuo_{i_1}\otdt\HEuo_{i_n}$, the action of $U$ is given by

\[
\zeta_1\otdt\zeta_n 
\mapsto\begin{cases}
\zeta_1\otimes(\zeta_2\otdt\zeta_{n-1})\otimes\zeta_n,&n\ge3,\,i_1=i_n=B \\
\zeta_1\otimes(\zeta_2\otdt\zeta_n)\otimes\oneh_B,&i_1=B,\,i_n=D \\
\oneh_B\otimes(\zeta_1\otdt\zeta_{n-1})\otimes\zeta_n,&i_1=D,\,i_n=B \\
\oneh_B\otimes(\zeta_1\otdt\zeta_n)\otimes\oneh_B,&i_1=i_n=D
\end{cases}
\]
in $L^2(B)\otimes\KEu\otimes L^2(B)$.
Since neither $B$ nor $D$ is equal to $\Cpx$ and  both are separable, the
Hilbert space $\KEu$ is separable and
infinite
dimensional.
The   equality~\eqref{eq:Ulam} is now easily obtained by verifying that 
\begin{equation}\label{eq.7.20}
U\lambda_N(b)=(\lambda_B(b)\;\oplus\;
\lambda_B(b)\otimes\id_{\KEu}\otimes\id_{L^2(B)})U.
\end{equation}
The equation~\eqref{eq:Urho} follows from equation~\eqref{eq:Ulam} and the fact
(see~\cite{V85}) that
for $\zeta_1\otdt\zeta_n\in\HEuo_{i_1}\otdt\HEuo_n$ as in~\eqref{eq:L2M},
\[
J_N(\zeta_1\otdt\zeta_n)=J_{i_n}\zeta_n\otdt J_{i_1}\zeta_1.
\]
\end{proof}

\begin{thm}\label{thm3.2}
Let $I$ be a finite or countable set containing at least two elements
and, for each $i\in I$, let $M_i$ be a diffuse separable von Neumann
algebra with
a normal faithful trace $\tau_i$. Let $A_i$ be a masa in $M_i$
and let $Q$ be a diffuse separable von Neumann algebra with
a normal faithful trace $\tau_Q$.
Let $\alpha_i\in(0,1)$ be such that $\sum_{i\in I}\alpha_i=1$.
Let
\[
M=\bigoplus_{i\in I}M_i
\]
with trace $\tau$ given by
\begin{equation}
\tau((x_i)_{i\in I})=\sum_{i\in I}\alpha_i\tau_i(x_i)
\end{equation}
and let
\begin{equation}\label{eq:A}
A=\bigoplus_{i\in I}A_i\subseteq M.
\end{equation}
Let
\begin{align}
(N_i,\tauhat_i)&=(M_i,\tau_i)*(Q,\tau_Q)\quad(i\in I) \label{eq:Ni} \\[1ex]
(N,\tauhat)&=(M,\tau)*(Q,\tau_Q) \label{eq:N}
\end{align}
be the free products of von Neumann algebras.
Then, for all $i\in I$, $N_i$ is a type ${\text{II}}_1$ factor and
$A_i$ is a masa in $N_i$.
Moreover $N$ is a type ${\text{II}}_1$ factor, $A$ is a masa in $N$ and
the Puk\'anszky invariants satisfy
\begin{equation}\label{eq:PukN}
\Puk_N(A)=\{\infty\}\cup\bigcup_{i\in I}\Puk_{N_i}(A_i).
\end{equation}
\end{thm}
\begin{proof}
We first note  that each $N_i$ and $N$ are factors 
in which respectively $A_i$ and $A$ are masas, by Theorem \ref{thm2.3} and
Remark \ref{rem2.6}.

Let $p_i\in A$ be the projection with entry $1$ in the $i$th component $A_i$
and entry $0$
in every other component of the direct sum~\eqref{eq:A}.
These are orthogonal, central projections in $M$.
For $i,j\in I$, let
\[
q_{ij}=\lambda_N(p_i)\rho_N(p_j)\in\Bc(L^2(N)).
\]
Then $q_{ij}$ is an element of $\Ac=(\lambda_N(A)\cup\rho_N(A))''$,
which is the center of $\Ac'$.
Since the strong operator sum $\sum_{i,j\in I}q_{ij}$ is equal to the identity,
it follows that the Puk\'anszky invariant $\Puk_N(A)$ is equal to the union
over
all $i$ and $j$ of those $n\in\Nats\cup\{\infty\}$ such that
$q_{ij}\Ac'(1-e_A)$ has
a nonzero part of type ${\text{I}}_n$.
Implicitly using the unitary of Lemma~\ref{lem3.1}
for the free product construction~\eqref{eq:N},
we identify $L^2(N)$ with
\[
L^2(M)\;\oplus\;L^2(M)\otimes\KEu\otimes L^2(M)
\]
and we make the identifications
\begin{align*}
\lambda_N(p_i)&=\lambda_M(p_i)\;\oplus\;
\lambda_M(p_i)\otimes\id_{\KEu}\otimes\id_{L^2(M)} \\
\rho_N(p_j)&=\rho_M(p_j)\;\oplus\;
\id_{L^2(M)}\otimes\id_{\KEu}\otimes\rho_M(p_j).
\end{align*}
If $i\ne j$, then $\lambda_M(p_i)\rho_M(p_j)=0$ and consequently $q_{ij}\Ac$
is identified with the algebra
\[
\big(\lambda_M(A_i)\otimes\id_{\KEu}\otimes\rho_M(A_j)\big)''
\subseteq\Bc(L^2(M_i)\otimes\KEu\otimes L^2(M_j)),
\]
which commutes with
$\id_{L^2(M_i)}\otimes\Bc(\KEu)\otimes\id_{L^2(M_j)}$. 
Therefore, $q_{ij}\Ac'$ is purely of type ${\text{I}}_\infty$.
Thus, $q_{ij}\Ac'$ contributes only $\infty$ to $\Puk_N(A)$.
On the other hand, in the case $i=j$, since $Ap_i=A_i$ and
$Mp_i=M_i$, we see that $q_{ii}\Ac$ is identified with the von Neumann algebra
generated by
\begin{multline*}
\{\lambda_{M_i}(a)\;\oplus\;
\lambda_{M_i}(a)\otimes\id_{\KEu}\otimes\id_{L^2(M_i)}\mid a\in A_i\}
\\
\cup
\{\rho_{M_i}(a)\;\oplus\;
\id_{L^2(M_i)}\otimes\id_{\KEu}\otimes\rho_{M_i}(a)\mid a\in A_i\}
\end{multline*}
in $\Bc(L^2(M_i)\oplus L^2(M_i)\otimes\KEu\otimes L^2(M_i))$.
Using the unitary $U$ from Lemma~\ref{lem3.1} in the case of the free
product~\eqref{eq:Ni},
we thereby identify $q_{ii}\Ac$ with the abelian von Neumann algebra $\Ac_i$
used to define $\Puk_{N_i}(A_i)$.
Thus, $q_{ii}\Ac'$ contributes exactly $\Puk_{N_i}(A_i)$ to $\Puk_N(A)$.
Taking the union over all $i$ and $j$ yields~\eqref{eq:PukN}.
\end{proof}

We now come to the main result of the paper. We let $\bb F_k$, $k\in
\{1,2,\ldots ,\infty\}$, denote the free group on $k$ generators, where $\bb
F_1$ is identified with $\bb Z$. To avoid discussion of separate cases below,
we adopt the convention that $\bb F_{k-r}$ means $\bb F_k$ when $k=\infty$
and $r\in \Nats$. 

\begin{thm}\label{thm3.3}
Let $S$ be an arbitrary subset of $\bb N$, let 
$k\in \{2,3, \ldots ,\infty\}$ be arbitrary, and let $\Gamma$ be an arbitrary
countable discrete group. Then there exists a strongly
singular masa $A\subseteq L(\bb F_k*\Gamma)$ whose Puk\'anszky invariant is
$S\cup
\{\infty\}$.
\end{thm}
\begin{proof} We first consider the case when $\Gamma$ is trivial. If $S$ is
empty, then we may take $A$ to be the masa
corresponding to one of the generators of $\bb F_k$.
This masa is strongly singular, \cite{SSgafa} and has Puk\'{a}nszky invariant
$\{\infty\}$, \cite{SSpuk}.
Thus we may assume that $S$ is nonempty. For each $n\geq 1$, let $G_n$ be the
I.C.C. group of Proposition \ref{prop2.5} with abelian subgroup $H_n$. The
subgroup
\[ \left\{\left(\begin{matrix} 1&x&y\\ 0&1&0\\ 0&0&1\end{matrix}\right)\colon 
\ 
x,y\in {\bb Q}\right\}
\]
is abelian and normal in $G_n$, and the quotient is isomorphic to $H_n$. As was
noted in \cite{SSpuk}, each $G_n$ is amenable, and $L(G_n)$ is the hyperfinite 
type ${\text{II}}_1$ factor $R$. If we write $M_n=L(G_n)\cong R$, then we have
the masa $A_n:=L(H_n)$ inside the hyperfinite type ${\text{II}}_1$ factor.
Then define $M_S = 
\bigoplus\limits_{i\in S} M_i$ and  $A_S = \bigoplus\limits_{i\in S} A_i$. We
see that $M_S$ is a direct sum of copies of $R$ and so is hyperfinite.
Consequently, $M_S*L(\bb F_1)$ is isomorphic to $L(\bb F_2)$, \cite{D}, while,
for $3\leq k \leq \infty$, we have the isomorphisms
\[M_S*L(\bb F_{k-1})\cong M_S*(L(\bb F_1)*L(\bb F_{k-2}))\cong (M_S*L(\bb
F_1))*L(\bb F_{k-2})\cong L(\bb F_2)*L(\bb F_{k-2})\cong L(\bb F_k).\]
The results of Section \ref{sect2} imply that $A_S$ is a strongly singular masa
in both $M_S$ and $M_S*L(\bb F_{k-1})$. From Theorem \ref{thm3.2},
\begin{equation}\label{eq3.20}
\Puk_{L(\bb F_k)}(A_S)=\{\infty\}\cup\bigcup_{i\in
S}\Puk_{M_i}(A_i)=S\cup\{\infty\},
\end{equation}
where the last equality comes from Proposition \ref{prop2.5}.

The case where the countable discrete group $\Gamma$ is included is essentially
the same. Simply observe that $M_S*(L(\bb F_{k-1})*L(\Gamma))\cong
(M_S*L(\bb F_{k-1}))*L(\Gamma)\cong L(\bb F_{k}))*L(\Gamma)$. 
\end{proof}

This theorem has an interesting parallel in \cite{NS}, where it was shown that
any subset of $\Nats\cup\{\infty\}$ which contains 1 can be $\Puk_R(A)$ for
some masa $A$ in the hyperfinite factor $R$.

When $M$ and $Q$ are type ${\text{II}}_1$ factors and $A$ is a masa in $M$ then
it is 
also a masa in $M*Q$ by Theorem \ref{thm2.3}.
In view of our results to this point, it is natural to ask whether 
\begin{equation}\label{eq3.21}
\Puk_{M*Q}(A)=\Puk_M(A)\cup\{\infty\}.
\end{equation}
  This is not true in general, as we now show. We
follow Proposition \ref{prop2.5}, but replace the groups there by
the  ones below which come from \cite[Example 5.1]{SSpuk}:
\begin{equation}\label{eq5.4}
G_n = \left\{\left(\begin{matrix} 1&x\\ 0&f\end{matrix}\right)\colon \ f\in
P_n, 
\quad x\in {\bb Q}\right\},\quad H_n = \left\{\left(\begin{matrix} 1&0\\ 
0&f\end{matrix}\right)\colon \ f\in P_n\right\}.
\end{equation}
The double cosets form a single equivalence class of $n$ elements and each
stabilizer subgroup is trivial (see \cite{SSpuk} for details). When the masa
$L(H_n)$ in $L(G_n)$
is viewed as a masa in $L(G_n*\Gamma)$ for a nontrivial countable discrete
group
$\Gamma$, the elements of $(G_n*\Gamma)\backslash G_n$ have trivial
stabilizer subgroups, and so we still have one equivalence class of double
cosets but now with infinitely many elements. We conclude that
$\Puk_{L(G_n)}(L(H_n))=\{n\}$, while
$\Puk_{L(G_n)*L(\Gamma)}(L(H_n))=\{\infty\}$. Thus \eqref{eq3.21} fails in
general.


\section{Masas in the hyperfinite factor}

\indent

In \cite{NS}, it was shown that any subset of $\Nats \cup\{\infty\}$ containing
1 could be the Puk\'anszky invariant of a masa in the hyperfinite type
${\rm{II}}_1$ factor $R$. Subsequently many other admissible subsets were
found in \cite{SSpuk}. Building on the examples of the latter paper, we  now
show that any subset of $\Nats \cup\{\infty\}$ containing $\{\infty\}$ is a
possible value of the invariant, exactly like the result obtained for the free
group factors in the previous section.
\begin{exm}\label{ex4.1}
We will need below the group $P_{\infty}$, defined by
\[P_{\infty} = \left\{\frac{p}q \colon \  p,q\in {\bb 
Z}_{\text{odd}}\right\},
\]
as in the proof of Proposition \ref{prop2.5}.

Let $S$ be an arbitrary subset of $\Nats$. We construct a strongly singular
masa $A$ in the the hyperfinite factor $R$ with $\Puk(A)=S\cup\{\infty\}$ as
follows. We will assume that $S$ is nonempty since $\{\infty\}$ is already a
known value (see \cite[Example 5.1]{SSpuk} with $n=\infty$). Let
$\{n_1,n_2,\ldots\}$ be a listing of the numbers in $S$, where each is
repeated infinitely often to ensure that the list is infinite. Then define a
matrix group $G$ by specifying the general group elements to be 
\begin{equation}\label{eq5.12}
\left(\begin{matrix} 1&x_1&x_2&x_3&\cdots\\ 0&f_12^{n_1k}&0&0&\cdots\\
0&0&f_22^{n_2k}&0&\cdots\\ 0&0&0& 
f_32^{n_3k}\\\vdots&\vdots&\vdots&&\ddots \end{matrix}\right),
\end{equation} 
where $k\in \Ints$, $x_j \in \bb Q$, $f_j\in P_{\infty}$, and the relations
$x_i \ne 0$ and $f_j\ne 1$ occur only finitely often. This makes $G$ a
countable group which is easily checked to be I.C.C. Moreover $G$ is amenable
since there is an abelian normal subgroup $N$ (those matrices with only 1's on
the diagonal) so that the quotient $G/N$ is isomorphic to the abelian subgroup
$H$ consisting of the diagonal matrices in $G$. Thus $R:=L(G)$ is the
hyperfinite factor and $A:=L(H)$ is a masa in $R$, as in the examples of
\cite[Section 5]{SSpuk}. For any finite nonempty subset $T$ of $\Nats$, let
$H_T$ be the subgroup of $H$ obtained by the requirements that $k=0$ and that
$f_i=1$ for $i\in T$.

Each nontrivial double coset is generated by a nontrivial element 
\begin{equation*}
x=\left(\begin{matrix} 1&x_1&x_2&x_3&\cdots\\ 0&1&0&0&\cdots\\
0&0&1&0&\cdots\\ 0&0&0&1\\\vdots&\vdots&\vdots&&\ddots \end{matrix}\right),
\end{equation*} 
of $N$, and we split into two
cases according to whether the number of nonzero $x_i$
is exactly 1 or is greater than 1. In the first case, suppose that $x_1$ is the
sole
nonzero value. We obtain $n_1$ distinct cosets generated by the elements
\begin{equation}\label{eq5.13}
\left(\begin{matrix} 1&2^k&0&0&\cdots\\ 0&1&0&0&\cdots\\ 0&0&1&0\\ 0&0&0& 
1\\\vdots&\vdots&&&\ddots \end{matrix}\right),\ \ 1\leq k\leq n_1,
\end{equation}
respectively,
whose stabilizer subgroups are all $\{(h,h^{-1}):\ h\in H_{\{1\}}\}$. A similar
result holds when the nonzero entry occurs in the $i^{{\rm{th}}}$ position:
$n_i$ distinct cosets with stabilizer subgroups $\{(h,h^{-1}):\ h\in
H_{\{i\}}\}$.

If $T$ is a finite subset of $\Nats$ with $|T|>1$, then two elements of $N$,
having nonzero entries respectively $x_i$ and $y_i$ for $i\in T$, generate the
same double coset precisely when there exist $k\in \Ints$ and $f_i\in
P_{\infty}$ such that $x_i=y_if_i2^{n_ik}$, $i\in T$. The stabilizer subgroup
in this case is $\{(h,h^{-1}):\ h\in H_{T}\}$. Thus the distinct stabilizer
subgroups are pairwise noncommensurable and \cite[Theorem 4.1]{SSpuk} allows
us to determine the Puk\'anszky invariant by counting the equivalence classes
of double cosets. In the first case we obtain the integers $n_i\in S$. In the
second case each equivalence class has infinitely many elements and so the
contribution is $\{\infty\}$, showing that
$
\Puk(A)=S\cup\{\infty\}$, as required.$\hfill\square$   

\end{exm}

\section{Representations of abelian C$^*$--algebras.}

We use the approach to direct integrals
found in Kadison and Ringrose~\cite[\S14.1]{KR},
because it fits well with our computations that
will follow in Section~\ref{sec:NScomp}.
Suppose $X$ is a compact Hausdorff space and $\pi:C(X)\to\Bc(\HEu)$
is a unital $*$--representation.
Then there is a Borel measure $\mu$ on $X$ so that $\HEu$ can be written as
a direct integral
\[
\HEu=\int^\oplus_X\HEu_x\,d\mu(x)
\]
and so that $\pi$ is a diagonal representation with, for all $f\in C(X)$ and
$x\in X$,
\[
\pi(f)_x=f(x)\id_{\HEu_x}.
\]
The multiplicity function $m(x)=\dim(\HEu_x)$ is $\mu$--measurable.
Of course such results are known in much greater generality.
The pair
$([\mu],m)$, where $m$ is taken up to redefinition on
sets of $\mu$--measure zero,
is a conjugacy invariant for $\pi$, and every such pair arises from
some unital $*$--representation of $C(X)$.
We will use Proposition~\ref{prop:di}
to find $\mu$ and $m$ in concrete examples.

\begin{defi}
We call $[\mu]$ the {\em measure class} of $\pi$ and $m$ the {\em multiplicity
function}
of $\pi$.
\end{defi}

We now consider what happens to the measure class and multiplicity function
under
certain natural constructions.

\begin{prop}\label{prop:ot}
For $i=1,2$ let $A_i=C(X_i)$ be an abelian, unital C$^*$--algebra and let
$\pi_i:A_i\to\Bc(\HEu_i)$ be a unital $*$--representation.
Let $[\mu_i]$ be the measure class and $m_i$ the multiplicity function of
$\pi_i$.
\begin{itemize}
\item[(i)]
Let $X$ be the disconnected union of $X_1$ and $X_2$, so that we identify
$A=C(X)$ with $A_1\oplus A_2$.
Let $\pi:A\to\Bc(\HEu_1\oplus\HEu_2)$ be the unital $*$--representation given
by
\[
\pi(f_1\oplus f_2)=\pi_1(f_1)\oplus\pi_2(f_2),\qquad(f_i\in A_i).
\]
Then the measure class of $\pi$ is $[\mu]$, where
\[
\mu(E_1\cup E_2)=\mu_1(E_1)+\mu_2(E_2),\qquad(E_i\subseteq X_i)
\]
and the multiplicity function $m$ of $\pi$ is such that the restriction
of $m$ to $X_i$ is $m_i$.
\item[(ii)]
Let $\pi=\pi_1\otimes\pi_2:A_1\otimes A_2\to\Bc(\HEu_1\otimes\HEu_2)$
be the tensor product representation of the tensor product C$^*$--algebra
$A_1\otimes A_2$, which we identify with $C(X_1\times X_2)$.
Then the measure class of $\pi$ is $[\mu_1\otimes\mu_2]$ and the multiplicity
function $m$ of $\pi$ is given by $m(x_1,x_2)=m_1(x_1)m_2(x_2)$.
\end{itemize}
\end{prop}
\begin{proof}
Part~(i) is obvious.

For part~(ii), by Lemma~14.1.23 of~\cite{KR}, we may without loss of generality assume
\begin{align*}
\mu_i&=\sum_{1\le n\le\infty}\mu_{i,n} \\
\HEu_i&=\bigoplus_{1\le n\le\infty}L^2(\mu_{i,n})\otimes\KEu_n, \\
\pi_i(f)&=\bigoplus_{1\le n\le\infty}M_f^{(i,n)}\otimes\id_{\KEu_n}\qquad(f\in A_i),
\end{align*}
for the family of mutually singular measures $(\mu_{i,n})_{1\le n\le\infty}$,
where $\KEu_n$ is a Hilbert space of dimension $n$ and
where $M_f^{(i,n)}$ is multiplication by $f$ on $L^2(\mu_{i,n})$.
Now the required formulas are easily proved.
\end{proof}

\begin{cor}\label{cor:otid}
Let $\pi:C(X)\to\Bc(\HEu)$ be a unital $*$--representation of an abelian
C$^*$--algebra
with measure class $[\mu]$ and multiplicity function $m$.
Let $\Vc$ be a Hilbert space with dimension $k$.
Then the representation $\pi\otimes\id:C(X)\to\Bc(\HEu\otimes\Vc)$ given by
$\pi\otimes\id(a)=\pi(a)\otimes\id_\Vc$ has measure class $[\mu]$ and
multiplicity
function $x\mapsto m(x)k$.
\end{cor}

\begin{prop}\label{prop:surj}
Suppose that $\phi:C(X)\to C(Y)$ is a surjective $*$--homomorphism of unital,
abelian
C$^*$--algebras.
We identify $Y$ with a closed subset of $X$ and $\phi$ with the restriction of
functions.
Let $\pi:C(Y)\to\Bc(\HEu)$ be a unital $*$--representation with measure class
$[\mu]$
and multiplicity function $m$.
Then the $*$--representation $\pit=\pi\circ\phi:C(X)\to\Bc(\HEu)$ has measure
class $[\mut]$
and multiplicity function $\mt$, where $\mut(E)=\mu(E\cap Y)$ and the
restriction of $\mt$
to $Y$ is $m$.
\end{prop}

Below, $\Hh$ will denote the dual of a locally compact abelian group $H$. We
recall the well known fact that $\Hh$ is compact when $H$ is discrete.

\begin{lem}
Let $H$ be a countable abelian group and let $X$ be a set on which $H$ acts
transitively.
Let $\pi$ be the $*$--representation of $C^*(H)\cong C(\Hh)$ on the Hilbert
space
$\ell^2(X)$ that results from this action.
Let $K\subseteq H$ be the stabilizer subgroup of any element
of $X$ under the action of $H$, and let
\[
K\oup=\{\gamma\in\Hh \mid\gamma(K)=\{1\}\}
\]
be the annihilator of $K$,
 a closed subgroup of the compact group $\Hh$.
Then the measure class of $\pi$ is supported on $K\oup$ and is equal
there to the class of Haar measure on $K\oup$.
The multiplicity function of $\pi$ is the constant function $1$.
\end{lem}
\begin{proof}
We identify $X$ with the quotient group $H/K$ equipped with the obvious action
of $H$.
Then $K\oup$ is the dual group of the quotient group $H/K$, so via the Fourier
transform yields $\ell^2(X)\cong L^2(K\oup,\lambda)$,
where $\lambda$ is Haar measure on $K\oup$.
This isomorphism intertwines the given representation of $C^*(H)=C(\Hh)$ on
$\ell^2(X)$
with the representation $\sigma$ of $C(\Hh)$ on $L^2(K\oup,\lambda)$ given by
\[
(\sigma(f)\xi)(\gamma)=f(\gamma)\xi(\gamma),
\qquad(f\in C(\Hh),\xi\in L^2(K\oup,\lambda),\gamma\in K\oup).
\]
\end{proof}

A special case of the above lemma is the following one.
\begin{lem}\label{lem:HaH}
Let $H$ be an abelian subgroup of a discrete I.C.C.\ group $G$
and let $\pi:C^*(H)\otimes C^*(H)\to\Bc(\ell^2(G))$ be the left--right
representation of the C$^*$--tensor product, given by
\[
\pi(\lambda_h\otimes\lambda_{h'})\delta_g=\delta_{hgh'}\qquad(h,h'\in H,\,g\in
G).
\]
We identify $C^*(H)\otimes C^*(H)$ with $C(\Hh\times\Hh)$.
Let $a\in G$ and let 
\[
K_a=\{(h_1,h_2)\in H\times H\mid h_1ah_2=a\}.
\]
Then the cyclic subrepresentation of $\pi$ on
$\ell^2(HaH)$ has measure class $[\mu]$,
where $\mu$ is concentrated on the annihilator subgroup
$(K_a)\oup\subset\Hh\times\Hh$
and is equal to Haar measure there.
The multiplicity function of $\pi$ is constantly equal to $1$.
\end{lem}

The next lemma shows how to write the direct sum of direct integrals of Hilbert
space
as a direct integral.

\begin{lem}\label{lem:Kdi}
Let $X$ be a $\sigma$--compact, locally compact Hausdorff space and let $\mu$
and $\mu'$ be
completions of $\sigma$--finite Borel measures on $X$.
Let separable Hilbert spaces $\HEu$ and $\HEu'$ be direct integrals of
$\{\HEu_p\}_{p\in X}$ and $\{\HEu'_p\}_{p\in X}$
over $(X,\mu)$ and $(X,\mu')$, respectively.
The Lebesgue decompositions yield
\begin{equation}\label{eq:Ldecomp}
\begin{alignedat}{2}
\mu&=\mu_0+\mu_1\quad&\text{with }&\mu_0\perp\mu',\quad\mu_1\ll\mu' \\
\mu'&=\mu'_0+\mu'_1\quad&\text{with }&\mu'_0\perp\mu,\quad\mu'_1\ll\mu.
\end{alignedat}
\end{equation}
Let $X=X_0\cup X_1=X'_0\cup X'_1$ be measurable partitions of $X$ such that
\begin{equation}\label{eq:X}
\mu_0(X_1)=\mu'(X_0)=\mu'_0(X'_1)=\mu(X'_0)=0.
\end{equation}
Thus, $\mu_0$ is concentrated on $X_0$, $\mu_1$ on $X_1$, $\mu'_0$ on $X'_0$
and $\mu'_1$ on $X'_1$.
Let $\nu=\mu+\mu'_0$, namely
\begin{equation}\label{eq:nu}
\nu(A)=\mu(A)+\mu'(A\cap X'_0).
\end{equation}
For $p\in X$ consider the Hilbert spaces
\[
\KEu_p=
\begin{cases}
0,&p\in X_0\cap X'_0 \\
\HEu'_p,&p\in X_1\cap X'_0 \\
\HEu_p,&p\in X_0\cap X'_1 \\
\HEu_p\oplus\HEu'_p,&p\in X_1\cap X'_1.
\end{cases}
\]
Then $\KEu:=\HEu\oplus\HEu'$ is the direct integral
of $\{\KEu_p\}_{p\in X}$ over $(X,\nu)$.  
Furthermore, if $a\in\Bc(\HEu)$ and $a'\in\Bc(\HEu')$
are decomposable with respect to the direct
integrals 
$\HEu=\int_X^\oplus\HEu_p\,d\mu(p)$ and $\HEu'=\int_X^\oplus\HEu'_p\,d\mu'(p)$,
respectively,
then $a\oplus a'\in\Bc(\HEu\oplus\HEu')$ is decomposable with respect to the
direct integral
\begin{equation}\label{eq:opdi}
\HEu\oplus\HEu'=\int_X^\oplus\KEu_p\,d\nu(p),
\end{equation}
with
\begin{equation}\label{eq:adecomp}
(a\oplus a')_p=
\begin{cases}
0,&p\in X_0\cap X'_0 \\
a'_p,&p\in X_1\cap X'_0 \\
a_p,&p\in X_0\cap X'_1 \\
a_p\oplus a'_p,&p\in X_1\cap X'_1.
\end{cases}
\end{equation}
\end{lem}
\begin{proof}
Let $f$ be the Radon--Nikodym derivative $d\mu'_1/d\mu$.
We may without loss of generality assume $f>0$ everywhere on $X_1\cap X'_1$.
Let $x\in\HEu$ and $x'\in\HEu'$ and set $\xt=x\oplus x'\in\KEu$.
For $p\in X$, we set
\[
\KEu_p\ni\xt(p)=
\begin{cases}
0,&p\in X_0\cap X'_0 \\  
x'(p),&p\in X_1\cap X'_0 \\ 
x(p),&p\in X_0\cap X'_1 \\  
x(p)\oplus f(p)^{1/2}x'(p),&p\in X_1\cap X'_1
\end{cases}
\]
A straightforward calculation shows
\[
\int_X\|\xt(p)\|^2\,d\nu(p)=
\|x\|^2+\| x'\|^2=\|\xt\|^2.
\]

Suppose that for all $p\in X$, $\ut(p)\in\KEu_p$ is such that
for all $\xt\in\KEu$, the function
$p\mapsto\langle\ut(p),\xt(p)\rangle$ is $\nu$--integrable,
and let us show there is $\yt\in\KEu$ such that $\yt(p)=\ut(p)$ for $\nu$-a.e.
$p\in X$.
For every $p\in X_1\cap X'_1$, let $Q_p:\KEu_p\to\HEu_p$ be the orthogonal
projection onto the first direct summand;
we denote by $I-Q_p:\KEu_p\to\HEu_p'$ the orthogonal projection onto
the second direct summand.
Set
\[
\HEu_p\ni u(p)=
\begin{cases}
0,&p\in X'_0 \\
\ut(p),&p\in X_0\cap X'_1 \\
Q_p\ut(p),&p\in X_1\cap X'_1.
\end{cases}
\]
Let $x\in\HEu$ and let $\xt=x\oplus0\in\KEu$.
Then for all $p\in X_1\cap X'_1$ we have
\[
\langle\ut(p),\xt(p)\rangle=
\langle\ut(p),Q_p\xt(p)\rangle=
\langle u(p),x(p)\rangle,
\]
while if $p\in X_0\cap X'_1$, then also
$\langle\ut(p),\xt(p)\rangle=\langle u(p),x(p)\rangle$.
Since the restrictions of $\mu$ and $\nu$ to $X'_1$ agree and since
$\mu(X'_0)=0$, it follows that the map $p\mapsto\langle u(p),x(p)\rangle$
is $\mu$--integrable.
Therefore, there is $y\in\HEu$ such that $y(p)=u(p)$ for $\mu$--a.e.\ $p\in X$.

Now let
\[
\HEu'_p\ni u'(p)=
\begin{cases}
0,&p\in X_0 \\
\ut(p),&p\in X_1\cap X'_0 \\
f(p)^{-1/2}(I-Q_p)\ut(p),&p\in X_1\cap X'_1.
\end{cases}
\]
Let $x'\in\HEu'$ and let $\xpt=0\oplus x'\in\KEu$,
so that
\[
\xpt(p)=
\begin{cases}
0,&p\in X_0 \\
x'(p),&p\in X_1\cap X'_0 \\
0\oplus f(p)^{1/2}x'(p),&p\in X_1\cap X'_1.
\end{cases}
\]
Then for all $p\in X_1\cap X'_1$ we have
\[
\langle\ut(p),\xpt(p)\rangle=
f(p)^{1/2}\langle(I-Q_p)\ut(p),x'(p)\rangle=
f(p)\langle u'(p),x'(p)\rangle,
\]
while for $p\in X_1\cap X'_0$ we have
$\langle\ut(p),\xpt(p)\rangle=\langle u'(p),x'(p)\rangle$.
Thus,
\begin{align*}
\int_X|\langle u'(p),x'(p)\rangle|\,d\mu'(p)
&=\int_{X_1\cap X'_1}|\langle u'(p),x'(p)\rangle|\,d\mu'(p)
+\int_{X_1\cap X'_0}|\langle u'(p),x'(p)\rangle|\,d\mu'(p) \\
&=\int_{X_1\cap X'_1}|\langle u'(p),x'(p)\rangle|f(p)d\mu(p)
+\int_{X_1\cap X'_0}|\langle u'(p),x'(p)\rangle|\,d\mu'(p) \\
&=\int_{X_1\cap X'_1}|\langle\ut(p),\xpt(p)\rangle|d\mu(p)
+\int_{X_1\cap X'_0}|\langle\ut(p),\xpt(p)\rangle|\,d\mu'(p) \\
&=\int_X|\langle\ut(p),\xpt(p)\rangle|\,d\nu(p)<\infty.
\end{align*}
Therefore, the function $p\mapsto\langle u'(p),x'(p)\rangle$ is
$\mu'$--integrable
and there is $y'\in\HEu'$ such that $u'(p)=y'(p)$ for $\mu'$-a.e. $p\in X$.
Let $\yt=y\oplus y'\in\KEu$.
Then $\yt(p)=\ut(p)$ for $\nu$-a.e.\ $p\in X$.
This proves the direct integral formula~\eqref{eq:opdi}.

The decomposability of $a\oplus a'$  with respect to~\eqref{eq:opdi}
and the validity of~\eqref{eq:adecomp} are now clear.
\end{proof}

A particular case of the above lemma is the following proposition,
which we will frequently use in computations.

\begin{prop}\label{prop:di}
Let $\pi:C(X)\to\Bc(\HEu)$ and $\pi':C(X)\to\Bc(\HEu')$ be unital
$*$--representations
whose measure classes are $[\mu]$ and $[\mu']$ and whose multiplicity functions
are
$m$ and $m'$, respectively.
Consider the Lebesgue decomposition as in~\eqref{eq:Ldecomp} and the partitions
so that we have~\eqref{eq:X}.
Then the representation $\pi\oplus\pi':C(X)\to\Bc(\HEu\oplus\HEu')$ has
measure class $[\nu]$ with $\nu=\mu+\mu'_0$ given by~\eqref{eq:nu},
and has multplicity function $\mt$ given by
\[
\mt(x)=
\begin{cases}
m'(x),&x\in X_1\cap X'_0 \\
m(x),&x\in X_0\cap X'_1 \\
m(x)+m'(x),&x\in X_1\cap X'_1.
\end{cases}
\]
\end{prop}

Before we state and prove the next result, we need to make some remarks which
will justify the calculations below. Consider an inclusion $M\subseteq N$ of
finite von~Neumann algebras where $N$ has a faithful finite normal trace
$\tau$. Then $L^1(N)$ denotes the completion of $N$ when equipped with the
norm $\|x\|_1=\tau(|x|)$, $x\in N$. The polar decomposition shows that
\begin{equation}\label{eq5.1000}
\|x\|_1=\sup \{|\tau(xy)|\colon y\in N,\ \|y\|\leq 1\},
\end{equation}
from which it follows that $|\tau(x)|\leq \|x\|_1$. Thus $\tau$ has a bounded
extension, also denoted $\tau$, to a linear functional on $L^1(N)$. The
$N$--bimodule properties of $N$ extend by continuity to $L^1(N)$ and the
relation $\tau(x\zeta)=\tau(\zeta x)$ for $x\in N$, $\zeta \in L^1(N)$,
follows by boundedness of $\tau$ on $L^1(N)$. Similarly, $\tau$ defines a
continuous linear functional on $L^2(N)$ by $\tau(\zeta)=\langle \zeta
,1\rangle$, for $\zeta \in L^2(N)$. If $\mathbb E$ is the unique trace
preserving conditional expectation of $N$ onto $M$ then it is also a
contraction when viewed as a map of $L^2(N)$ onto $L^2(M)$. If $x\in N$, then
\begin{align}
 \|\mathbb E(x)\|_1 &=\sup\{|\tau(\mathbb E(x)y)|\colon y\in M,\ \|y\|\leq 1\}
=\sup\{|\tau(\mathbb E(xy))|\colon y\in M,\ \|y\|\leq 1\}\nonumber\\ 
&=\sup\{|\tau(xy)|\colon y\in M,\ \|y\|\leq 1\}\leq \|x\|_1,\label{eq5.1001}
\end{align}
using the module properties of $\mathbb E$. Thus $\mathbb E$ has a bounded
extension to a contraction of $L^1(N)$ to $L^1(M)$. The module property 
$ \mathbb E(m_1\zeta m_2)=m_1\mathbb E(\zeta )m_2$, for $\zeta \in L^1(N)$ and
$m_1, m_2 \in M$, follows by $\|\cdot\|_1$--continuity of $\mathbb E$. The
bilinear map $\Psi \colon N\times N \to N$, defined by $\Psi(x,y)=xy$,
satisfies
\begin{align}
\|\Psi(x,y)\|_1&=\sup \{|\tau(xyz)|\colon z\in N,\ \|z\|\leq 1\}\nonumber\\
&=\sup \{|\langle zx,y^* \rangle |\colon z\in N,\ \|z\|\leq 1\}\nonumber\\
&\leq \|x\|_2\|y\|_2,
\end{align}
and so $\Psi$ extends to a jointly continuous map, also denoted $\Psi$, from
$L^2(N)\times L^2(N)$ to $L^1(N)$. Since $\Psi$ is the product map at the
level of $N$, we will write $\zeta \eta$ for $\Psi(\zeta, \eta)\in L^1(N)$
when $\zeta, \eta\in L^2(N)$. Moreover, the adjoint on $N$ extends to an
isometric conjugate linear map $\zeta \mapsto \zeta^*$ on both $L^1(N)$ and
$L^2(N)$, agreeing with $J$ in the latter case. Any relation that holds on $N$
will extend by continuity to the appropriate $L^p(N)$, where $p=1{\text{ or
}}2$. For example, we have
\begin{equation}
\langle \zeta, \eta \rangle =\tau (\eta^* \zeta)=\tau (\zeta \eta^*),\ \ \
\zeta, \eta \in L^2(N).
\end{equation}
We will apply these remarks with $N=\overline{B}$ and $M=\overline{A}$ below.

\begin{prop}\label{prop:lamnu}
Let $A=C(X)$ be embedded as a unital C$^*$--subalgebra
of a separable C$^*$--algebra $B$ and let $\tau$ be a faithful,
tracial state on $B$.
Let $\nu$ be the Borel measure on $X$ such that
\[
\tau(a)=\int_Xa(x)\,d\nu(x),\qquad(a\in A).
\]
Let $\lambda$ denote the representation of $A$ on $L^2(B,\tau)$
by left multiplication.
Then the measure class of $\lambda$ is $[\nu]$.
\end{prop}
\begin{proof}
Let $[\sigma]$ denote the measure class of $\lambda$.
Let $\lambda_A$ denote the left action of $A$ on $L^2(A,\tau)$.
Now $\lambda_A$ is a direct summand of $\lambda$ and the measure
class of $\lambda_A$ is easily seen to be $[\nu]$.
By Proposition~\ref{prop:di}, $\nu\ll\sigma$.
In order to show $\sigma\ll\nu$, it will suffice to show that
whenever $\pi$ is a cyclic subrepresentation of $\lambda$, then
the measure class of $\pi$ is absolutely continuous with respect to $\nu$.
Indeed, if $\sigma\not\ll\nu$, then letting $X_0\subset X$ be a
set of positive measure such that the restriction of $\sigma$ to $X_0$ is
singular to $\nu$, there is $\zeta\in L^2(B,\tau)$ such that
$\int_{X_0}\|\zeta(x)\|^2\,d\nu(x)>0$.

Let $\zeta\in L^2(B,\tau)$.
Let $\Bbar$ denote the s.o.--closure of $B$
acting via the Gelfand--Naimark--Segal representation on $L^2(B,\tau)$,
and let $\Abar$ denote the s.o.--closure of $A$ acting via $\lambda$ on
$L^2(B,\tau)$.
Let $\mathbb E:\Bbar\to\Abar$ denote the $\tau$--preserving conditional
expectation
and also the extension
\[
\mathbb E:L^1(\Bbar,\tau)\to L^1(\Abar,\tau),
\]
which exists by the preceding remarks.
The measure class of $\pi$ is $[\rho]$, where for all $a\in C(X)$,
$\langle a\zeta,\zeta\rangle=\int_X a(x)\,d\rho(x)$.
Assume $a\ge0$.
Then, using the discussion before this proposition, we have
\[
\langle a\zeta,\zeta\rangle=\tau(\zeta^*a\zeta)=\tau(a\zeta\zeta^*)
=\tau(\mathbb E(a\zeta\zeta^*))=\tau(a\mathbb E(\zeta\zeta^*)).
\]
Since $0\le \mathbb E(\zeta\zeta^*)\in L^1(A,\tau)\cong L^1(\nu)$,
we have
\[
\langle a\zeta,\zeta\rangle=\int_Xa(x)f(x)\,d\nu(x)
\]
for some $f\in L^1(\nu)$, $f\ge0$.
Consequently, $d\rho=f\,d\nu$ and $\rho$ is absolutely continuous with respect
to $\nu$.
\end{proof}

\begin{prop}\label{prop:fp}
For $i=1,2$, let $B_i$ be a unital, separable C$^*$--algebra having faithful,
tracial states $\tau_i$ and with $\dim(B_i)\ge2$.
Let $A=C(X)$ be unitally embedded as a C$^*$--subalgebra of $B_1$.
Let $\nu$ be the measure on $X$ such that
\[
\tau_1(a)=\int_Xa(x)\,d\nu(x),\qquad(a\in C(X)).
\]
Let $\lambda_1,\rho_1:A\to\Bc(L^2(B_1))$ be the left and right
actions of $A$ on $L^2(B_1):=L^2(B_1,\tau_1)$.
Let $\pi_1:A\otimes A\to\Bc(L^2(B_1))$ be the $*$--representation of the
C$^*$--tensor product $A\otimes A$
given by $\pi_1(a_1\otimes a_2)=\lambda_1(a_1)\rho_1(a_2)$.
Let $(B,\tau)=(B_1,\tau_1)*(B_2,\tau_2)$ be the reduced free product of
C$^*$--algebras.
Let $\lambda,\rho:A\to\Bc(L^2(B))$ be the left and right actions of $A$ on
$L^2(B):=L^2(B,\tau)$ and let $\pi:A\otimes A\to\Bc(L^2(B))$ be given by
$\pi(a_1\otimes a_2)=\lambda(a_1)\rho(a_2)$.
Then $L^2(B_1)\subseteq L^2(B)$ is a reducing subspace for $\pi(A\otimes A)$.
Let $\pit$ be the representation of $A\otimes A$ on $L^2(B)\ominus L^2(B_1)$
obtained from $\pi$ by restriction, so that
\[
\pi=\pi_1\oplus\pit.
\]
We identify $A\otimes A$ with $C(X\times X)$.
Then the measure class of $\pit$ is $[\nu\otimes\nu]$ and the multiplicity
function of $\pit$ is constantly $\infty$.
\end{prop}
\begin{proof}
This proof is at bottom quite similar to the proof of Lemma~\ref{lem3.1},
and just as in that case, we begin by decomposing the free product Hilbert
space.
From the construction of the reduced free product C$^*$--algebra, we have
\[
L^2(B)=\Cpx\oneh\oplus
\bigoplus_{\substack{n\in\Nats \\ i_1,\ldots,i_n\in\{1,2\} \\ i_j\ne i_{j+1}}}
\HEu\oup_{i_1}\otdt\HEu_{i_n}\oup,
\]
where $\HEu\oup_i=L^2(B_i)\ominus\Cpx\oneh$.
Therefore,
\begin{equation}\label{eq:LBom}
L^2(B)\ominus L^2(B_1)=
\bigoplus_{k=0}^\infty\bigg(\begin{aligned}[t]
&\HEu\oup_2\otimes(\HEu\oup_1\otimes\HEu\oup_2)^{\otimes k}\quad\oplus\quad
\HEu\oup_2\otimes(\HEu\oup_1\otimes\HEu\oup_2)^{\otimes k}\otimes\HEu\oup_1 \\
&\quad\oplus\quad
 \HEu\oup_1\otimes\HEu\oup_2\otimes(\HEu\oup_1\otimes\HEu\oup_2)^{\otimes k} \\
&\quad\oplus\quad
\HEu\oup_1\otimes
  \HEu\oup_2\otimes(\HEu\oup_1\otimes\HEu\oup_2)^{\otimes
k}\otimes\HEu\oup_1\bigg).
\end{aligned}
\end{equation}
Writing $L^2(B_1)=\HEu\oup_1\oplus\Cpx$, for each $k$ we make the
obvious identification of the direct sum of the
four direct summands on the right hand side
of~\eqref{eq:LBom} with
\[
L^2(B_1)\otimes\HEu\oup_2\otimes(\HEu\oup_1\otimes\HEu\oup_2)^{\otimes
k}\otimes L^2(B_1).
\]
Let
\[
\KEu=\bigoplus_{k=0}^\infty\HEu\oup_2\otimes(\HEu\oup_1\otimes\HEu\oup_2)^{\otimes
k}.
\]
Then we have the unitary
\[
U:L^2(B)\ominus L^2(B_1)\to L^2(B_1)\otimes\KEu\otimes L^2(B_1)
\]
that gives
\[
U(\pit(a_1\otimes a_2))U^*=\lambda_1(a_1)\otimes\id_\KEu\otimes\rho_1(a_2),
\qquad(a_1,a_2\in A).
\]
By Corollary~\ref{cor:otid}, it will suffice to show that the representation
\[
\lambda_1\otimes\rho_1:A\otimes A\to\Bc(L^2(B_1)\otimes L^2(B_1))
\]
has measure class $[\nu\otimes\nu]$.
By Proposition~\ref{prop:lamnu}, $\lambda_1$ has measure class $\nu$.
Since $\rho_1$ is unitarily equivalent to $\lambda_1$ it also has measure class
$\nu$,
so by Proposition~\ref{prop:ot},
$\lambda_1\otimes\rho_1$ has measure classs $[\nu\otimes\nu]$.
\end{proof}


\section{Computations of  invariants}
\label{sec:NScomp}

\begin{para}\label{para:NS}
Neshveyev and St\o rmer~\cite{NS} considered the conjugacy invariant for
a masa $A$ in a II$_1$ factor $M$ derived from writing a direct integral
decomposition of the left--right action,
\begin{equation}\label{eq:repab}
a\otimes b\mapsto aJb^*J,
\end{equation}
where $J$ is the anti--unitary conjugation on $L^2(M)$ given by
$J\ah=(a^*)\hat{\;}$,
of the C$^*$--tensor
product $A\otimes A$ on $L^2(M)$.
We will now review this invariant.
Choosing a separable and weakly dense C$^*$--subalgebra $\Afr=C(Y)$ of $A$,
we may write $A=L^\infty(Y,\nu)$ for a compact Hausdorff
space $Y$ and a completion of a Borel measure $\nu$ on $Y$.
Let $\pi:\Afr\otimes\Afr\to\Bc(L^2(M))$ denote the restriction
of the left--right action~\eqref{eq:repab} to the
C$^*$--tensor product $\Afr\otimes\Afr$, which we identify with $C(Y\times Y)$
in the usual way.
Let $[\eta]$ be the measure class and $m$ the multiplicity function of $\pi$.
We will always take $\eta$ to be a finite measure.
Then $[\eta]$ is invariant
under the flip $(x,y)\mapsto(y,x)$ of $Y\times Y$ and
the projection of $[\eta]$ onto the first and second coordinates is $[\nu]$.
Neshveyev and St\o{}rmer~\cite{NS} observed that $(Y,[\eta],m)$ is a conjugacy
invariant
of $A\subseteq M$, in the sense that if $A\subseteq M$ and $B\subseteq N$ are
masas
and if there is an isomorphism $M\to N$ taking $A$ onto $B$, then
(for any choice of separable, weakly dense C$^*$--subalgebras of $A$ and $B$),
there
is a transformation of measure spaces, $F:(Y_A,[\nu_A])\to(Y_B,[\nu_B])$
such that $(F\times F)_*([\eta_A])=[\eta_B]$ and $m_B\circ(F\times F)=m_A$
($\eta_A$--almost everywhere).
We will refer to the equivalence class of $(Y,[\eta],m)$ under such
transformations
as the {\it {\MM}} invariant of the masa $A\subseteq M$.
In fact, Neshveyev and St\o{}rmer showed more, namely that
the equivalence class of $(Y,[\eta],m)$ is a complete invariant for
the pair $(A,J)$ acting on $L^2(M)$.
They also showed that the Puk\'anszky invariant of $A\subseteq M$ is
precisely the set of essential values
of the multiplicity function $m$ taken on the complement
of the diagonal $\Delta(Y)$ in $Y\times Y$.
\end{para}

\begin{exm}\label{ex:A}
Let $n\in\Nats\cup\{\infty\}$ and consider the matrix groups
\[
G_n=\left\{
\left(\begin{matrix}
f&x \\ 0&1  
\end{matrix}\right)
\biggm|
f\in P_n,\,x\in\bb{Q}
\right\},\qquad
H_n=\left\{
\left(\begin{matrix}
f&0 \\ 0&1  
\end{matrix}\right)
\biggm|
f\in P_n
\right\}\subseteq G_n,
\]
where we have the subgroups of the multiplicative group of nonzero rational
numbers
\[
P_\infty=\{\frac pq\mid p,q\in\Ints^*,\,p,q\text{ odd}\}
\]
and, for $n$ finite,
\[
P_n=\{f2^{kn}\mid f\in P_\infty,\,k\in\Ints\}.
\]
Then, as shown in~\cite[Ex.\ 5.1]{SSpuk},
$L(G_n)$ is the hyperfinite II$_1$ factor and $L(H_n)$ is a strongly singular
masa
in $L(G_n)$ with Puk\'anszky invariant $\{n\}$.
Moreover, the \MM invariant of $L(H_n)\subseteq L(G_n)$
is the equivalence class of $(\Hh_n,[\mu_n],m)$,
where $\mu_n$ is the sum of Haar measure on $\Hh_n\times\Hh_n$
and Haar measure on the diagonal subgroup $\Delta(\Hh_n)$,
and where the multiplicity function $m$ takes value $1$ on $\Delta(\Hh_n)$
and  $n$ on its complement.
\end{exm}
\begin{proof}
This follows from the double decomposition of $G_n$ as double cosets over $H_n$
(see~\cite[Ex.\ 5.1]{SSpuk}), Lemma~\ref{lem:HaH} and
Proposition~\ref{prop:di}.
\end{proof}

\begin{exm}\label{ex:B}
Let $n\in\Nats$.
With $H_n\subseteq G_n$ and $H_\infty\subseteq G_\infty$ as in
Example~\ref{ex:A},
$L(H_n\times H_\infty)$ is a strongly singular masa in $L(G_n\times G_\infty)$
whose \MM invariant is the equivalence class of
\[
(\Hh_n\times\Hh_\infty,[\eta],m),
\]
where $\eta$ is the sum of
\begin{itemize}
\item[(i)] Haar measure on $\Hh_n\times\Hh_\infty\times\Hh_n\times\Hh_\infty$
\item[(ii)] Haar measure on the subgroup
\begin{equation}\label{eq:Dn}
D_n=\{(\alpha,\beta_1,\alpha,\beta_2)
 \mid\alpha\in\Hh_n,\,\beta_1,\beta_2\in\Hh_\infty\}
\end{equation}
\item[(iii)] Haar measure on the subgroup
\[
D_\infty=\{(\alpha_1,\beta,\alpha_2,\beta)
 \mid\alpha_1,\alpha_2\in\Hh_n,\,\beta\in\Hh_\infty\}
\]
\item[(iv)] Haar measure on the diagonal subgroup
$\Delta(\Hh_n\times\Hh_\infty)$
\end{itemize}
and where the multiplicity function $m$ is given by
\[
m(\gamma)=
\begin{cases}
1,&\gamma\in\Delta(\Hh_n\times\Hh_\infty) \\
n,&\gamma\in D_n\backslash\Delta(\Hh_n\times\Hh_\infty) \\
\infty,&\text{else.}
\end{cases}
\]
\end{exm}
\begin{proof}
This follows from Example~\ref{ex:A} and Proposition~\ref{prop:ot}.
\end{proof}

\begin{exm}\label{ex:C}
Let $n\in\Nats$ and let $\Gamma$ be any nontrivial finite or countably infinite
group.
Let
\[
H_n\times H_\infty\subseteq G_n\times G_\infty\subset(G_n\times
G_\infty)*\Gamma,
\]
with $H_n\times H_\infty\subseteq G_n\times G_\infty$ as in Example~\ref{ex:B}
above.
Then $L(H_n\times H_\infty)$ is a singular masa in $L((G_n\times
G_\infty)*\Gamma)$,
whose 
\MM invariant is the equivalence class of
\[
(\Hh_n\times\Hh_\infty,[\eta],m),
\]
where $\eta$ and $m$ are exactly as in Example~\ref{ex:B}.
\end{exm}
\begin{proof}
This follows from Example~\ref{ex:B} and applications of
Propositions~\ref{prop:fp} and~\ref{prop:di}.
\end{proof}

If $\Gamma$ is taken to be infinite amenable, then $L((G_n\times
G_\infty)*\Gamma)$
is isomorphic to the free group factor $L(\Fb_2)$, by~\cite{D93a} (see
also~\cite{D}).

\begin{exm}\label{ex:D}
Let
\[
A=L(H_\infty)\oplus L(H_n\times H_\infty)\subseteq N=L(G_\infty)\oplus
L(G_n\times G_\infty).
\]
Consider the normal, faithful, tracial state
\[
\tau_N(x_1\oplus x_2)=\frac12(\tau_{G_\infty}(x_1)+\tau_{G_n\times
G_\infty}(x_2))
\]
on $N$.
Let $Q$ be any diffuse von Neumann algebra with separable predual and a normal
faithful state $\tau_Q$ and let
\[
(M,\tau_M)=(N,\tau_N)*(Q,\tau_Q)
\]
be the free product.
By Theorem~\ref{thm2.3}, $M$ is a II$_1$ factor and $A\subseteq M$ is a
strongly
singular
masa.
The \MM invariant of $A\subseteq M$ 
is the equivalence class of
$(X,[\sigma],m)$,
where $X$ is the disconnected sum of $X_1=\Hh_\infty$ and
$X_2=\Hh_n\times\Hh_\infty$,
where $\sigma$ is the sum of the measures
\begin{itemize}
\item[(i)] $\mu_\infty$, as described in Example~\ref{ex:A},
supported on $X_1\times X_1\subseteq X\times X$
\item[(ii)] $\eta$, as described in Example~\ref{ex:B},
supported on $X_2\times X_2\subseteq X\times X$
\item[(iii)] $\nu\otimes\nu$ on $X\times X$, where the measure $\nu$ on $X$ is
the sum of Haar measure on the (dual) group $X_1$ and Haar measure on $X_2$,
\end{itemize}
and where $m$ takes the value $1$ on the diagonal $\Delta(X)$,
the value $n$ on $D_n\subseteq X_2\times X_2\subseteq X\times X$,
with $D_n$ as given in equation~\eqref{eq:Dn}, and is equal to $\infty$
elsewhere.
\end{exm}
\begin{proof}
Let $\Afr=C^*(H_\infty)\oplus C^*(H_n\times H_\infty)\subseteq A$.
We find the measure class and multiplicity function of the left--right
representation
of the C$^*$--algebra $\Afr\otimes\Afr$ on $L^2(N,\tau_N)$,
by using Propositions~\ref{prop:ot}
and~\ref{prop:surj}.
Then we find the measure class of the left--right representation of
$\Afr\otimes\Afr$
on $L^2(M,\tau_M)$ by using Propositions~\ref{prop:fp} and~\ref{prop:di}.
\end{proof}

If $Q$ is taken to be the hyperfinite II$_1$ factor, then $M$ in
Example~\ref{ex:D}
is isomorphic to the free group factor $L(\Fb_2)$, by~\cite{D}.
Thus, Examples~\ref{ex:C} and~\ref{ex:D} provide two constructions of
masas
in the free group factor $L(\Fb_2)$,
both having Puk\'anszky invariant $\{n,\infty\}$.
We will distinquish these two masas
using the \MM invariant, or actually a formally weaker
invariant
derived from it.

\begin{para}\label{para:NSoffdiag}
Let $(Y,[\eta],m)$ arise as in the definition of the \MM
invariant
of a masa $A\subseteq M$.
As already mentioned, $m$ takes the value $1$ on the diagonal $\Delta(Y)$,
and one easily sees that the restriction of $\eta$ to $\Delta(Y)$ is equivalent
to the measure $\nu$ as in~\ref{para:NS}, when $\Delta(Y)$ is identified with
$Y$
in the obvious way.
Therefore, the restrictions of $m$ and $\eta$ to the complement of $\Delta(Y)$
contain the same information as $(Y,[\eta],m)$.
\end{para}

\begin{lem}\label{lem:suppproj}
Let $Q$ be a von Neumann algebra having normal faithful traces $\tau_1$ and
$\tau_2$
and let $A\subseteq Q$ be a von Neumann subalgebra.
Let $E_i:Q\to A$ denote the $\tau_i$--preserving conditional expectation onto
$A$,
($i=1,2$).
If $x\in Q$ and $x\ge0$, then the support projections of $E_1(x)$ and $E_2(x)$
agree.
\end{lem}
\begin{proof}
Let $p_i\in A$ be such that the support projection of $E_i(x)$ is $1-p_i$.
Then $0=p_1E_1(x)p_1=E_1(p_1xp_1)$, so $p_1xp_1=0$.
But then $p_1E_2(x)p_1=E_2(p_1xp_1)=0$, so $p_1\le p_2$.
By symmetry, $p_2\le p_1$.
\end{proof}

\begin{para}\label{para:Bool}
Let $M$ be a II$_1$ factor and $A\subseteq M$ a masa.
Choose a triple $(Y,[\eta],m)$ belonging to the \MM invariant
of $A\subseteq M$, with $\eta$ finite, as considered in~\ref{para:NS}.
Recall that $A=L^\infty(Y,\nu)$ is embedded in $L^\infty(Y\times Y,\eta)$
as functions constant in the second coordinate.
Let $E:L^\infty(Y\times Y,\eta)\to A$ denote the conditional expectation
that preserves integration with respect to $\eta$.
Given $n\in\Nats\cup\{\infty\}$, let $P_n\in L^\infty(Y\times Y,\eta)$ be the
characteristic
function of the set where the multiplicity function $m$ takes the value $n$ off
of the
diagonal $\Delta(Y)$,
and let $q_n(A)=q_n(A,M)$
be the support projection of the conditional expectation $E(P_n)$.
By Lemma~\ref{lem:suppproj}, $q_n(A)$ is independent of the choice of $\eta$ in
the measure
class $[\eta]$.
Moreover, if $F:(Y_A,[\nu_A])\to(Y_B,[\nu_B])$ is the transformation of measure
spaces
considered in~\ref{para:NS}, then we have $q_n(A,M)=q_n(B,N)\circ F$ for
the corresponding support projections.
Therefore, if $A\subset M$ and $B\subset N$ are masas that are conjugate by and
isomorphism
from $M$ to $N$, then it induces an isomorphism from $A$ to $B$ that sends
$q_n(A)$ to $q_n(B)$ for 
all $n$.
\end{para}

Fix $n\in\Nats$.
In Example~\ref{ex:C}, take $\Gamma=\Ints$, so that
$M=L((G_n\times G_\infty)*\Gamma)=L(\Fb_2)$
and let $A_{\ref{ex:C}}=L(H_n\times H_\infty)$ be the masa of $L(\Fb_2)$
obtained there.
In Example~\ref{ex:D}, take $Q$ to be the hyperfinite II$_1$ factor so that
$M=L(\Fb_2)$
and let $A_{\ref{ex:D}}$ be the masa of $L(\Fb_2)$ obtained there.
\begin{thm}
The masas $A_{\ref{ex:C}}$ and $A_{\ref{ex:D}}$ in $L(\Fb_2)$
both have Puk\'anszky invariant $\{n,\infty\}$,
but are non-conjugate.
\end{thm}
\begin{proof}
The values of the Puk\'anszky invariant can be read off from the
\MM invariants, which were computed in Examples~\ref{ex:C}
and~\ref{ex:D}.
The derived invariant $q_n$ from~\ref{para:Bool} above in these cases becomes
$q_n(A_{\ref{ex:C}})=1$ while $q_n(A_{\ref{ex:D}})=1_{X_2}$, the characteristic
function
of $X_2\subseteq X$, which is not the identity of $A_{\ref{ex:D}}$.
\end{proof}

\end{document}